\documentclass[11pt]{article}

\usepackage{amsmath, amssymb, amsthm, calc}  

\title{Theorem on best Diophantine approximations for linear forms.
                              \thanks{The research was supported by
                              RFBR (grant $\textup N^{\circ}$ 06--01--00518).
                              The first author was also supported by
                              the grant of the President of Russian Federation
                              $\textup N^\circ$ MK--4466.2008.1.
                              The Russian version of this paper
                              is submitted to the Proceedings of the
                              Department of Mathematics and Mechanics,
                              Moscow State University.
                             }}
\author{Oleg N. German, Nikolay G. Moshchevitin \\ \\ \\}
\date{\begin{flushright} \it
        Dedicated to Professor Viktor Antonovich Sadovnichiy, \\
        member of RAS, \\
        in occasion of his 70-th birthday.
      \end{flushright}}

\theoremstyle{definition}

\theoremstyle{remark}

\theoremstyle{plain}
\newtheorem{theorem}{Theorem}
\newtheorem{lemma}{Lemma}

\newtheorem{proposition}{Proposition}
\newtheorem{corollary}{Corollary}

\DeclareMathOperator{\conv}{conv}

\renewcommand{\geq}{\geqslant}
\renewcommand{\leq}{\leqslant}

\newcommand{\R}{\mathbb{R}}
\newcommand{\Z}{\mathbb{Z}}
\newcommand{\Q}{\mathbb{Q}}

\begin{document}

\maketitle

\begin{abstract}
We prove a new quantitative result on the degeneracy of the dimension of the subspace spanned by
the best Diophantine approximations for a linear form.
\end{abstract}

\section{Best approximations.}

Let $\alpha_1,\dots,\alpha_r$ be real numbers. Suppose that $1,\alpha_1,\dots,\alpha_r$ are
linearly independent over the rationals. For an integer point $m=(m_0,m_1,\dots,m_r)\in\Z^r$ we
define
\[ \zeta(m)=m_0+m_1\alpha_1+\dots+m_r\alpha_r. \]
A point $m=(m_0,m_1,\dots,m_r)\in\Z^{r+1}\setminus\{0\}$ is defined to be \emph{a best
approximation (in the sense of linear form)} (briefly \emph{b.a.}) if
\begin{equation} \label{eq:b_a_definition}
  \zeta(m)=\min_n\|\zeta(n)\|,
\end{equation}
where $\|\cdot\|$ denotes the distance to the nearest integer and the minimum is taken over all the
integer vectors $n=(n_0,n_1,\ldots,n_r)\in\Z^r$ such that
\[ 0<\max_{1\leq j\leq r}|n_j|\leq\max_{1\leq j\leq r}|m_j|. \]
All the best approximations form a sequence of points $m_\nu=(m_{0,\nu},m_{1,\nu},\dots,m_{r,\nu})$
with increasing $\max_{1\leq j\leq n}|m_{j,\nu}|$. It should be noticed that sometimes the points
$-m_\nu$ are also called best approximations. From this point of view, as can be seen from
\eqref{eq:b_a_definition}, in each pair $\pm m_\nu$ we chose as $m_\nu$ the point with positive
$\zeta(m_\nu)$.

Let us denote
\begin{equation} \label{0}
  \zeta_\nu=\zeta(m_\nu),\quad M_\nu=\max_{1\leq j\leq n}|m_{j,\nu}|.
\end{equation}
Then
\[ \zeta_1>\zeta_2>\cdots>\zeta_\nu>\zeta_{\nu+1}>\cdots \]
and
\[ M_1<M_2<\cdots<M_\nu<M_{\nu +1}<\cdots. \]
It follows from the Minkowski convex body theorem that $\zeta_\nu M_{\nu+1}^r\leq1$. Define
$\Delta_\nu^r$ to be the determinant of the matrix formed by the coefficients of $r+1$ consecutive
best approximations:
\[ \Delta_\nu^r=\left|
   \begin{array} {cccc}
     m_{0,\nu } & m_{1,\nu } &  \dots & m_{r,\nu } \cr
     \dots & \dots & \dots & \dots \cr
     m_{0,\nu +r} & m_{1,\nu +r} & \dots & m_{r,\nu +r}
   \end{array}\right|. \]

In the case $r=1$ one can easily see from the theory of continued fractions that
$\Delta_\nu^1=(-1)^{\nu-1}$ for every $\nu$.

H.\,Davenport and W.\,M.\,Schmidt were the first to prove that in the case $n=2$ there are
infinitely many values of $\nu$, for which $\Delta_\nu^2\neq0$ (see Lemma 3 from \cite{DS}). This
result follows from the Minkowski convex body theorem.

N.\,G.\,Moshchevitin \cite{MDAN},\cite{ME} obtained the following result. Given $r\geq3$, there is
an uncountable set of $r$-tuples $(\alpha_1,\dots,\alpha_r)$ of real numbers, linearly independent
with the unit over the rationals, such that for each $r$-tuple in this set all the points $m_\nu$
of the corresponding sequence of best approximations, starting with some $\nu$, lie in a certain
three-dimensional sublattice $\Lambda(\alpha_1,\ldots,\alpha_r)$ of the lattice $\Z^{r+1}.$ In this
paper we prove a more precise version of this result:

\begin{theorem} \label{t:zeta_convergence_implies_degeneracy}
Suppose that points
\[ m_\nu=(m_{0,\nu},m_{1,\nu},\dots,m_{r,\nu})\in\Z^{r+1},\quad\nu=1,2,3,\ldots \]
form the sequence of b.a. for $\alpha_1,\dots,\alpha_r$. Take an integer $k\geq1$ and consider the
integer points
\begin{equation} \label{sequ}
m_\nu^*=(m_{0,\nu},m_{1,\nu},\dots,m_{r,\nu},\underbrace{0,\ldots,0}_{k\text{ раз}})\in
\Z^{r+k+1},\quad\nu=1,2,3,\ldots
\end{equation}
Suppose that the series
\begin{equation} \label{er}
\sum_{\nu=1}^\infty M_{\nu+1}^{r+k}(\log M_{\nu+1})^{\delta_k}\zeta_\nu,\quad\delta_k=
\begin{cases}
1,\quad k =1 \cr
0,\quad k \geq 2
\end{cases},
\end{equation}
(where $M_\nu$  and $\zeta_\nu$ are defined by \eqref{0}) converges. Then for almost all
$(\beta_1,\ldots,\beta_k)\in\R^k$ (in the sense of Lebesgue measure) the sequence of best
approximations to the $(r+k)$-tuple $(\alpha_1,\ldots,\alpha_r,\beta_1,\ldots,\beta_k)$ differs
from the sequence $m_\nu^*$ by at most a finite set of integer points.
\end{theorem}

It is shown in the next section that in the case $r\geq2$ the series \eqref{er} may converge,
indeed. Thus it follows from Theorem \ref{t:zeta_convergence_implies_degeneracy} with $r=2$ that
there is a $(k+2)$-tuple $(\alpha_1,\alpha_2,\beta_1,\ldots,\beta_k)$ and an integer $\nu_0$ such
that all the b.a. $m_\nu^*,\nu\geq\nu_0$ lie in a three-dimensional subspace. The proof of Theorem
\ref{t:zeta_convergence_implies_degeneracy} is close to the proof of Lemma 3 from \cite{DS1}. We
prove Theorem \ref{t:zeta_convergence_implies_degeneracy} in Section \ref{sec:theorem_proof}.

\section{Khintchine's $\psi$-singular linear forms.}

Suppose that $\psi(y)=o(y^{-r})$, $y\to+\infty$. An $r$-tuple $(\alpha_1,\dots,\alpha_r)$ is said
to be $\psi$-singular (in the sense of linear form) if for every $T>1$ the Diophantine inequalities
\[ \|m_1\alpha_1+\dots+m_r\alpha _r\|<\psi(T),\quad0<\max_{1\leq j\leq r}|m_j|\leq T \]
have a solution in integer $r$-tuple $m=(m_1,\ldots,m_r)$. For $r\geq2$ A.~Khintchine \cite{HINS}
(see also Ch. 5, \S 7 of \cite{Cassil}) proved the existence of $\psi$-singular $r$-tuples for an
arbitrary function $\psi$.

It is easy to verify that an $r$-tuple $(\alpha_1,\dots,\alpha_r)$ of real numbers, linearly
independent with the unit over $\Q$,
is $\psi $-singular if and only if for all positive integer $\nu $ the following inequality is
valid
\begin{equation} \label{1}
  \zeta _\nu  \leq \psi ( M_{\nu +1}).
\end{equation}
Theorem \ref{t:zeta_convergence_implies_degeneracy} leads immediately to the following result.

\begin{theorem} \label{t:psi_convergence_implies_degeneracy}
Let $r\geq 2$ and let $\alpha_1,\dots,\alpha_r$ form a $\psi$-singular $r$-tuple. Suppose that the
series
\begin{equation} \label{er1}
 \sum_{\nu=1}^{\infty}M_\nu^{r+k}(\log M_\nu)^{\delta_k}\psi(M_\nu)
\end{equation}
converges. Then for almost all $(\beta_1,\ldots,\beta_k)\in\R^k$ the sequence of all the best
approximations to the $(r+k)$-tuple $(\alpha_1,\ldots,\alpha_r,\beta_1,\ldots,\beta_k)$ differs
from the sequence \eqref{sequ} by at most a finite set of integer points.
\end{theorem}

Now we discuss the convergence of the series \eqref{er}.

\begin{lemma} \label{l:exponential_growth}
For any $r$-tuple $(\alpha_1,\ldots,\alpha_r)$ of real numbers, linearly independent with the unit
over $\Q$, and for any positive integer $\nu$ one has
\begin{equation} \label{bestest}
  M_{\nu+2^{2r+1}-2^{r+1}}\geq2M_\nu.
\end{equation}
\end{lemma}

We prove Lemma \ref{l:exponential_growth} in Section \ref{sec:lemma_proof}. The proof is close to
that of a similar statement for the simultaneous approximations (see Theorem 2.2 from \cite{LAG} or
Lemma 1 from \cite{MOSHEV}). It is based on the pigeon hole principle. Lemma
\ref{l:exponential_growth} shows that the coefficients of the best approximations increase
exponentially and hence for any positive $\varepsilon$ the series
\[ \sum_{\nu =1}^{\infty} \frac{1}{(\log M_\nu)^{1+\varepsilon}} \]
converges. So we have the following

\begin{corollary} \label{cor:log_psi_implies_degeneracy}
Let $r\geq2$ and
\[ \psi(y)=\frac{1}{y^{r+k}(\log y)^{\delta_k+1+\varepsilon}} \]
with some positive $\varepsilon$. Let an $r$-tuple $(\alpha_1,\dots,\alpha_r)$ be $\psi$-singular.
Then for almost all $(\beta_1,\ldots,\beta_k)\in\R^k$ the sequence of best approximations for the
$(r+k)$-tuple $(\alpha_1,\ldots,\alpha_r,\beta_1,\ldots,\beta_k)$ differs from the sequence
\eqref{sequ} by at most a finite set of integer points.
\end{corollary}

One may ask if the inequality \eqref{bestest} can be essentially improved for $\psi$-singular
$r$-tuples expecting a-priori that for example under certain conditions on $\psi$ the series
\begin{equation}
\sum_{\nu =1}^{\infty} \frac{1}{\log M_\nu} \label{ser}
\end{equation}
should converge in the case of a $\psi$-singular $r$-tuple $(\alpha_1,\ldots,\alpha_r)$. This is
however not the case, for the methods of the paper \cite{ME} allow to prove the following

\begin{proposition}
Suppose that $r\geq 2$. Then, given an arbitrary function $\psi(y)=o(y^{-r})$, there is a
$\psi$-singular $r$-tuple $(\alpha_1,\ldots,\alpha_r)$ such that the series \eqref{ser} diverges.
\end{proposition}

But those same methods allow to prove that there are some special $\psi$-singular $r$-tuples for
which it is possible to improve Corollary \ref{cor:log_psi_implies_degeneracy}.

\begin{proposition} \label{prop:psi_singularity_plus_linear_independence}
Suppose that $r\geq 2$. Then, given an arbitrary function $\psi(y)=o(y^{-r})$, there is a
$\psi$-singular $r$-tuple $(\alpha_1,\ldots,\alpha_r)$ such that for every $\nu$ the vectors
$m_\nu,m_{\nu +1},\ldots,m_{\nu+r}$ are linearly independent.
\end{proposition}

\begin{corollary} \label{cor:loglog_psi_implies_degeneracy}
Suppose that
\begin{equation} \label{psi}
\psi (y) = \frac{1}{ y^{r+k} (\log y)^{\delta_k}(\log\log y)^{1+\varepsilon}}
\end{equation}
with some  positive $\varepsilon$. Suppose also that a $\psi$-singular $r$-tuple
$(\alpha_1,\ldots,\alpha_r)$ satisfies the condition that every $\nu$ vectors
$m_\nu,m_{\nu+1},\ldots,m_{\nu+r}$ are linearly independent. Then for almost all
$(\beta_1,\ldots,\beta_k)\in\R^k$ the sequence of best approximations for the $(r+k)$-tuple
$(\alpha_1,\ldots,\alpha_r,\beta_1,\ldots,\beta_k)$ differs from the sequence \eqref{sequ} by at
most a finite set of integer points.
\end{corollary}

\begin{proof}
Observe that under the conditions of Proposition
\ref{prop:psi_singularity_plus_linear_independence} we have the inequality
\begin{equation} \label{beste}
\zeta_\nu\geq\frac{1}{(r+1)!M_{\nu+r}^r}
\end{equation}
for every $\nu$. Indeed, consider the parallelepiped
\[ \Pi=\Big\{(x_0,x_1,\ldots,x_r)\in\R^{r+1}\ \Big|\, \max_{1\leq j\leq r}|x_j|\leq M_{\nu+r},\
   |x_0+x_1\alpha_1+\cdots+x_r\alpha_r|\leq\zeta_\nu\Big\}. \]
Observe that
\[ m_\nu,m_{\nu+1},\ldots,m_{\nu+r}\in\Pi. \]
The convex hull
\[ \mathcal O=\conv(\pm m_\nu,\pm m_{\nu+1},\ldots,\pm m_{\nu+r})\subset\Pi \]
is an integer $(r+1)$-dimensional polytope. For its $(r+1)$-dimensional measure we have the
following upper bound:
\begin{equation} \label{upper}
  \mu\mathcal O\leq\mu\Pi=2^{r+1}\zeta_\nu M_{\nu+r}^r.
\end{equation}
Since $\mathcal O$ is a lattice polytope,
\begin{equation} \label{lower}
  \mu\mathcal O\geq\frac{2^{r+1}}{(r+1)!}\,.
\end{equation}
Combining \eqref{upper} and \eqref{lower} we get \eqref{beste}.

Now \eqref{beste} together with \eqref{1} and the special choice of $\psi$ by \eqref{psi} leads to
the estimate
\[ \frac{1}{(r+1)!M_{\nu+r}^r}\leq\zeta_\nu\leq\psi(M_{\nu+1})=
   \frac{1}{M_{\nu+1}^{r+k}(\log M_{\nu+1})^{\delta_k}(\log\log M_{\nu+1})^{1+\varepsilon}}\,. \]
So, for $\nu$ large enough ($\nu\geq\nu_0=\nu_0(k,r)$) we have
\[ M_{\nu+r}\geq M_{\nu+1}^{1+\frac{k}{r}}. \]
Hence
\[ \log\log M_{\nu}\geq c\nu \]
with some positive $c=c(k,\alpha_1,\ldots,\alpha_r)$, which implies the convergence of the series
\eqref{er1}. The statement of the Corollary now follows from Theorem
\ref{t:psi_convergence_implies_degeneracy}.
\end{proof}

\section{Proof of Theorem \ref{t:zeta_convergence_implies_degeneracy}.} \label{sec:theorem_proof}

It is sufficient to prove that for almost every $(\beta_1,\ldots,\beta_k) \in [0,1]^k$ there is a
$\nu_0$ such that for all $\nu\geq \nu_0$ one has
\begin{equation} \label{eq:min_geq_zeta}
  \min|m_0+m_1\alpha_1+\cdots+m_r\alpha_r+m_{r+1}\beta_1+\cdots+m_{r+k}\beta_k|\geq\zeta_\nu
\end{equation}
where the minimum is taken over all the integer points $(m_0,m_1,\ldots,m_{r+k})$ such that
\[ \max_{1\leq j\leq r+k}|m_j|\leq M_{\nu+1},\quad|m_{r+1}|+\cdots+|m_{r+k}|\neq0. \]
This condition \eqref{eq:min_geq_zeta} is equivalent to the fact that for all integers
$m_0,\ldots,m_{r+k}$ such that
\begin{equation} \label{p}
  \max_{1\leq j\leq r}|m_j|\leq M_{\nu+1},\quad0<\max_{r+1\leq j\leq r+k}|m_j|\leq M_{\nu+1}
\end{equation}
one has \begin{equation} \label{3}
  m_{r+1}\beta_1+\cdots+m_{r+k}\beta_k \not\in  J_\nu(m_0,m_1,\ldots,m_r)
\end{equation}
where
\[ J_\nu(m_0,m_1,\ldots,m_r) = (-m_0 -m_1\alpha_1-\cdots-m_r\alpha_r- \zeta_\nu, - m_0
   -m_1\alpha_1-\cdots-m_r\alpha_r+ \zeta_\nu) \]
is an interval of length $2\zeta_\nu$. The condition \eqref{3} in its turn means that the  distance
between the point $(\beta_1,\ldots,\beta_k)\in [0,1]^k$ and the  subspace
\[ \Big\{ (x_1,\ldots,x_k)\in\R^k \,\Big|\
   m_{r+1}x_1+\cdots+m_{r+k}x_k=-m_0-m_1\alpha_1-\cdots-m_r\alpha_r \Big\} \]
is not less than $\zeta_\nu(m_{r+1}^2+\cdots+m_{r+k}^2)^{-1/2}$. Put
\[ \Omega_\nu (m_0,\ldots,m_{r+k})=\Big\{ (x_1,\ldots,x_k)\in[0,1]^k \,\Big|\
   m_{r+1}x_1+\cdots+m_{r+k}x_k\in J_\nu(m_0,m_1,\ldots,m_r) \Big\}. \]
We thus must prove that for almost every $(\beta_1,\ldots,\beta_k)\in[0,1]^k$ there is a $\nu_0$
such that
\[ (\beta_1,\ldots,\beta_k)\not\in\bigcup_{\nu\geq\nu_0}
   \left( \bigcup_m \Omega_\nu(m_0,\ldots,m_{r+k})\right) \]
where the inner union is taken over all the integers $m_0,m_1,\ldots,m_{r+k}$ satisfying the
condition \eqref{p} and the inequality  $|m_0|\leq (r+k+1)M_{\nu+1}$.

For the measure of $\Omega_\nu(m_0,\ldots,m_{r+k})$ we have
\[ \mu  \Omega_\nu (m_0,\ldots,m_{r+k})\leq
   \frac{2(k+r+1)k^{k/2}\zeta_\nu}{\sqrt{m_{r+1}^2+\cdots+m_{r+k}^2}}\,. \]
Hence
\[ \begin{split}
     \mu\left(\bigcup_m\Omega_\nu(m_0,\ldots,m_{r+k})\right)\leq
     \qquad\qquad\quad\ & \\
     \leq2(k+r+1)k^{k/2}\zeta_{\nu} \times (2M_{\nu+1}+1)^{r+1} \times &
     \sum_{0<\max_{1\leq j\leq k}|m_{r+j}|\leq M_{\nu+1}}
     (m_{r+1}^2+\cdots+m_{r+k}^2)^{-1/2}\leq \\
     & \qquad\qquad\qquad\qquad\qquad\quad\,
     \leq cM_{\nu+1}^{r+k}(\log M_{\nu+1})^{\delta_k}\zeta_{\nu}
   \end{split} \]
with some positive $c=c(k,\alpha_1,\ldots,\alpha_r)$. The series \eqref{er} converges, so it
remains to apply the Borel-Cantelli lemma.

\section{Proof of Lemma \ref{l:exponential_growth}.} \label{sec:lemma_proof}

Suppose that \eqref{bestest} is not true for some $\nu$. Then all the best approximations $m_j$,
$\nu\leq j\leq\nu+2^{2r+1}-2^{r+1}$ lie in the set
\[ \Pi=\Big\{ (x_0,x_1,\ldots,x_r)\in\R^{r+1}\ \Big|\ M_\nu\leq\max_{1\leq j\leq r}|x_j|<2M_{\nu},
   \ |x_0+x_1\alpha_1+\cdots+x_r\alpha_r|\leq\zeta_\nu  \Big\}. \]
This set can be covered by $2^{2r+1}-2^{r+1}$ half-open parallelepipeds of the form
\[ \begin{split}
     \Big\{ (x_0,x_1,\ldots,x_r)\in\R^{r+1}\ \Big|\ M_{\nu}/2<x_j-x_j^*\leq M_{\nu}/2, & \ \
     j=1,\ldots,r, \\
     & \eta\cdot(x_0+x_1\alpha_1+\cdots+x_r\alpha_r)\in[0,\zeta_\nu] \Big\}
   \end{split} \]
with some $(x_0^*,\ldots,x_r^*)\in\R^{r+1}$ and $\eta\in\{-1,+1\}$. By the pigeon hole principle
one of these parallelepipeds contains at least two distinct points $m_i$ and $m_j$. Hence
\[ m_j-m_i\in\Big\{ (x_0,x_1,\ldots,x_r)\in \R^{r+1}\ \Big|\ \max_{1\leq j\leq r}|x_j|<M_{\nu},\
  |x_0+x_1\alpha_1+\cdots+x_r\alpha_r|\leq \zeta_\nu \Big\}, \]
which contradicts the fact that $m_\nu$ is a best approximation and proves the Lemma.

\vskip1.5cm

\noindent
Oleg N. {\sc German} \\
Moscow Lomonosov State University \\
Vorobiovy Gory, GSP--1 \\
119991 Moscow, RUSSIA \\
\emph{E-mail}: {\fontfamily{cmtt}\selectfont german@mech.math.msu.su, german.oleg@gmail.com}

\vskip1.0cm

\noindent Nikolay  G. {\sc Moshchevitin} \\
Moscow Lomonosov State University \\
Vorobiovy Gory, GSP--1 \\
119991 Moscow, RUSSIA \\
\emph{E-mail}: {\fontfamily{cmtt}\selectfont moshchevitin@mech.math.msu.su,
moshchevitin@rambler.ru}

\end{document}